\documentclass[11pt]{article}
\usepackage[margin=4.46cm]{geometry}
\usepackage{graphicx}
\usepackage{epsfig}
\usepackage{amsmath}
\usepackage{amssymb}
\usepackage{caption}

\newcommand{\commentout}[1]{}

\def \Rset {{\mathbb R}}

\def \Zset {{\mathbb Z}}
\def \Qset {{\mathbb Q}}
\def \Nset {{\mathbb N}}

\newcommand{\nit}{\noindent}

\newcommand{\be}{\begin{equation}}
\newcommand{\ee}{\end{equation}}
\newcommand{\ba}{\begin{eqnarray}}
\newcommand{\ea}{\end{eqnarray}}
\newcommand{\bi}{\begin{itemize}}
\newcommand{\ei}{\end{itemize}}
\newcommand{\br}{\begin{eqnarray}}
\newcommand{\er}{\end{eqnarray}}

\newcommand{\qed}{\mbox{$\square$}\newline}

\newtheorem{theo}{Theorem}[section]
\newtheorem{defin}{Definition}[section]

\newtheorem{lem}{Lemma}[section]
\newtheorem{cor}{Corollary}[section]
\newtheorem{rmk}{Remark}[section]

\begin{document}
\title{{\Large\bf{High Degeneracy of Effective Hamiltonian in Two Dimensions}} }
\author{ Yifeng Yu\thanks{
This work is partially supported by NSF grant DMS-2000191.
}}

\date{}
\maketitle

\begin{abstract}
Consider the effective Hamiltonian $\overline H(p)$ associated with the mechanical Hamiltonian $H(p,x)={1\over 2}|p|^2+V(x)$. We prove that for generic $V$, $\overline H$ is piecewise 1d in a dense open set in two dimensions using Aubry-Mather theory. 
\end{abstract}

\section{Introduction}

Assume that $H=H(p,x)\in C(\Rset^n\times \Rset^n)$ is $\Zset^n$-periodic in $x$ and uniformly coercive in $p$, i.e
$$
\lim_{|p|\to +\infty}\min_{x\in \Rset^n}H(p,x)=+\infty.
$$
For each $\epsilon>0$, let $u^{\epsilon}\in C(\mathbb{R}^n\times [0,\infty))$ be the viscosity solution to the following Hamilton-Jacobi equation
\be\label{HJ-ep}
\begin{cases}
u_{t}^{\epsilon}+ H\left(Du^{\epsilon}, {x\over \epsilon}\right)=0 \quad &\text{in $\mathbb{R}^n\times (0,\infty)$},\\
u^{\epsilon}(x,0)=g(x) \quad &\text{on $\mathbb{R}^n$}.
\end{cases}
\ee
It was proved by Lions, Papanicolaou and Varadhan \cite{LPV} that $u^{\epsilon}$, as $\epsilon\to 0$, converges locally uniformly to $u$, the solution of the effective equation,
\be\label{HJ-lim}
\begin{cases}
u_t+{\overline H}(Du)=0 \quad &\text{in $\mathbb{R}^n\times (0,\infty)$},\\
u(x,0)=g(x) \quad &\text{on $\mathbb{R}^n$}. 
\end{cases}
\ee
Here $\overline H:\Rset^n\to \Rset$ is the so called {\it ``effective Hamiltonian}", which is determined by the following cell problem: for any $p\in \Rset^n$, there exists a unique number $\overline H$ such that the equation
\be\label{cell}
H(p+Dv,x)=\overline H(p)
\ee
has a periodic viscosity solution. When $H$ is convex in $p$, the effective Hamiltonian is convex and has a variational formulation
\be\label{inf-max}
\overline H(p)=\inf_{\phi\in C^1(\Bbb T^n)}\max_{x\in \Bbb T^n}H(p+D\phi(x),x).
\ee

Although there is a lot of literature regarding homogenization of Hamilton-Jacobi equations in various settings, not much is known about finer properties of the effective Hamiltonian $\overline H$ due to lack of tools. In this paper, we focus on the mechanical Hamiltonian 
$$
H(p,x)={1\over 2}|p|^2+V(x).
$$
Here $V$ is assumed to be $C^k(\Rset^n)$ for $k\geq 2$ and $\Zset ^n$-periodic. The following properties hold in any dimension.
\medskip

\nit{\bf $\bullet$ Property 1: Quadratic growth.}
$$
{1\over 2}|p|^2+\min_{\Rset^n} V\leq \overline H(p)\leq {1\over 2}|p|^2+\max_{\Rset^n} V.
$$
\nit {\bf $\bullet$ Property 2: Minimum Value}
$$
\min_{\Rset^n}\overline H=\max_{\Rset^n}V.
$$
More interestingly, for quite general $V$, the minimum level set
$$
F_0=\{p\in \Rset^n|\ \overline H=\max_{\Rset^n}V\}
$$
is an $n$-dimensional convex set (\cite{Co}).

\medskip

\nit {\bf $\bullet$ Property 3: Strict convexity along non-tangential direction} Using techniques from weak KAM theory, it was proved in \cite{EG} that $\overline H$ is not linear along any direction that is not tangent to its level set. In particular, this implies that if 
$$
\overline H\left({p_1+p_2\over 2}\right)={1\over 2}\overline H(p_1)+{1\over 2}\overline H(p_2),
$$
then 
$$
H(\lambda p_1+(1-\lambda)p_2)=\overline H(p_1) \quad \text{for all $\lambda\in [0,1]$}.
$$

Although $\overline H$ inherits some global features of ${1\over 2}|p|^2$, its local properties could be drastically different from those of ${1\over 2}|p|^2$. In this paper, we will prove that for generic $V$, $\overline H$ is piecewise 1d on a dense open set using Aubry-Mather theory. Precisely speaking,

\begin{theo}\label{main2} For any $k\geq 2$, there is a residual subset $\mathcal{G}$ of $C^k(\Bbb T^2)$ such that for every $V\in \mathcal{G}$, there exists a sequence of bounded open sets $\{O_i\}_{i\geq 1}$ in $\Rset^2$ such that

(1)
$$
O_V=\cup_{i=1}^{\infty}O_i
$$
is a dense open set in $\Rset^2$;

(2) For each $i\in \Nset$, there exist a unit vector $q_i\in \Rset^2$ and a convex function $f_i:\Rset\to \Rset$ such that 
$$
\overline H_V(p)=f_i(q_i\cdot p) \quad \text{in $O_i$}.
$$
Here $\overline H_V$ is the effective Hamiltonian associated with ${1\over 2}|p|^2+V$.
\end{theo}
In general, the above conclusion might not be true. For example, consider the separable case when $V(x)=h(x_1)+g(x_2)$. In dynamical system, some generic perturbation mechanisms based on Baire property in general topology have often been employed to filter out those exceptional situations. However, it is usually impossible to tell whether a concrete example is generic or not. 

$\mathcal{G}$ in the above theorem is the intersection of a sequence of dense open sets of $C^k(\Bbb T^2)$. There are two choices of $\mathcal{G}$. One choice is to directly use the residual set in Corollary 1.2 of \cite{BC} whose existence is established under certain abstract frameworks of convex analysis. The other choice is the $\mathcal{G}$ constructed in the appendix (\ref{residue}), which is weaker in the sense of dynamical system but is more explicit and enough for our purpose. Moreover, the result is expected to hold for more general Hamiltonian. In this paper, for clarity of presentation, we will only focus on the mechanical Hamiltonian that is interesting enough.

\medskip

\nit{\bf Notation and terminology:} 

\medskip

\nit $\bullet$ $\Bbb T^n=\Rset^n/ \Zset^n$ represents the $n$-dimensional flat torus. $C^{k}(\Bbb T^n)$ is the set of all $\Zset^n$-periodic $C^k(\Rset^n)$ functions. 

\medskip

\nit $\bullet$ A vector $q\in \Rset^n$ is called a {\it rational vector} if there exists $\lambda\in \Rset$ such that $\lambda q\in \Zset^n$. 

\medskip

\nit $\bullet$ $(m,n)\in \Zset^2$ is called {\it irreducible} if $|m|$ and $|n|$ are relatively prime. 

\medskip

\nit $\bullet$ A curve $\xi:\Rset\to \Bbb T^n$ is called periodic if there exists $T>0$ such that
$$
\xi(t+T)=\xi(t) \quad \text{for all $t\in \Rset$}.
$$
$T$ is called a period. If $T_0>0$ is the minimal period of $\xi$ and $\xi$ is lifted to $\Rset^n$, then 
$$
\xi(T_0)-\xi(0)=(m,n)\in \Zset^2
$$
is the first homology class of $\xi$. 

\medskip

\nit $\bullet$ Denote by $L(q,x):\Rset^n\times \Rset^n\to \Rset$ the Lagrangian
$$
L(q,x)=\sup_{p\in \Rset^n}\{q\cdot p-H(p,x)\}.
$$

\medskip

\nit $\bullet$ For $p_1$, $p_2\in \Rset^n$, let
$$
[p_1,p_2]=\{tp_1+(1-t)p_2|\ t\in [0,1]\}
$$
be the line segment connecting $p_1$ and $p_2$.

\section{Preliminary} For readers' convenience, in this section, we give a brief review about some basic knowledge and relevant results of Aubry-Mather theory and the weak KAM theory. See \cite{B1}, \cite{W-E}, \cite{EG} and \cite{F} for more details. Our presentation will be mainly from PDE point of view. Many parts are very close to the standard theory of Hamilton-Jacobi equations. The major difference lies on the classical Aubry-Mather theory, which is based on 2d topology and can not be captured by PDE approaches.

Let $\Bbb T^n=\Rset^n/\Zset^n$ be the $n$-dimensional flat torus and $H(p,x)\in C^{2}(\Rset^n\times \Rset^n)$ be a Hamiltonian satisfying
\begin{itemize}
\item[(H1)] (Periodicity) $x\mapsto H(p,x)$ is $\Zset^n$-periodic;

\item[(H2)] (Uniform convexity) There exists $\theta>0$ such that for all $\eta=(\eta_1,...,\eta_n)\in \Rset ^n$,
and $(p,x)\in \Rset^n \times \Rset^n$,
$$
\sum_{i,j=1}^{n}\eta_{i}{\partial ^2H\over \partial p_i\partial
p_j}\eta_{j}\geq \theta|\eta|^2.
$$
\end{itemize}
A major goal in dynamical system is to understand long time behaviors of trajectories of the Hamiltonian system
$$
\begin{cases}
\dot x(t)=D_pH(p,x)\\
\dot p(t)=-D_xH(p,x).
\end{cases}
$$
When the Hamiltonian $H$ is a small perturbation of the integrable case, the famous KAM theory based on analytic approaches says that most trajectories lie on invariant tori and hence integrable. For general $H$, when $n=2$, the classical Aubry-Mather theory provides nice description of structures of action minimizing trajectories based on topological approaches. See section 2.3 for more details. In \cite{Ma1}, Mather has extended the Aubry-Mather theory to higher dimensions through variational methods. The weak KAM theory reveals interesting connections between Mather's theory and solutions to the cell problem (\ref{cell}). In fact, in the language of PDE, the classical KAM theory can be formulated as: in the perturbation case, for ``most" $Q=D\overline H(p)$, the cell problem has a unique smooth solution $v$ up to a constant and the corresponding invariant torus is given by
$$
\Bbb T_{Q}= \{(q,x)\in \Rset ^n\times \Rset^n:\ p+Dv(x)=D_{q}L(q,x)\}. 
$$

\subsection{Aubry set and Man\'e set}

Let $v$ be a solution to the cell problem (\ref{cell}) . We say that a curve $\gamma:\Rset\to \Rset^n$ is a {\it global charaterstics} associated with $v$ if for all $t_1<t_2$ and $u(x)=p\cdot x+v$
$$
\int_{t_1}^{t_2}L(\dot \gamma,\gamma)+\overline H(p)\,ds=u(\gamma(t_2))-u(\gamma(t_1)).
$$
Write
$$
J_v=\cup_{\gamma} \{(\dot \gamma(t), \gamma(t))|\ t\in \Rset, \ \text{$\gamma$ is a global characteristics of $v$}\}.
$$
Such a $\gamma$ is also called a {\it $(v,L,\overline H(p))$-calibrated curve} in \cite{F}. According to classical theory in Hamilton-Jacobi equations \cite{Lions}, that $\gamma$ is a global characteristics of $v$ is equailvant to saying that $v$ is differentiable along $\gamma$ and for all $t\in \Rset$
\be\label{diff-char}
p+Dv=Du(\gamma(t))=D_qL(\dot \gamma(t), \gamma(t)).
\ee
Moreover, owing to Lemma \ref{calib}, every global characteristics is an absolutely minimizing curve with respect to $L(q,x)+\overline H(p)$. So it satisfies the Euler–Lagrange equation
 equation
\be\label{E-L}
{d(D_qL(\dot \gamma(t), \gamma(t)))\over dt}=D_xL(\dot \gamma(t), \gamma(t)).
\ee
A curve $\gamma: \Rset\to \Rset^n$ is called an {\it absolutely minimizing curve with respect to
$L(q,x)+c$} if for any
$-\infty<s_2<s_1<\infty$, $-\infty<t_2<t_1<\infty$ and $\eta\in AC([s_1,s_2], \Rset^n)$ subject to $\eta (s_2)=\gamma (t_2)$ and $\eta (s_1)=\gamma (t_1)$
the following inequality holds, 
\be\label{abs}
\int_{s_1}^{s_2} \left( L(\dot \eta(s),\eta(s))+c\right)\,ds\geq 
\int_{t_1}^{t_2} \left( L(\dot\gamma(s),\gamma (s))+c\right)\,ds.
\ee
Here $ AC([a,b], S)$ stands for the set of absolutely continuous curves $[a,b]\to S$. 
\begin{center}
\includegraphics[scale=0.5]{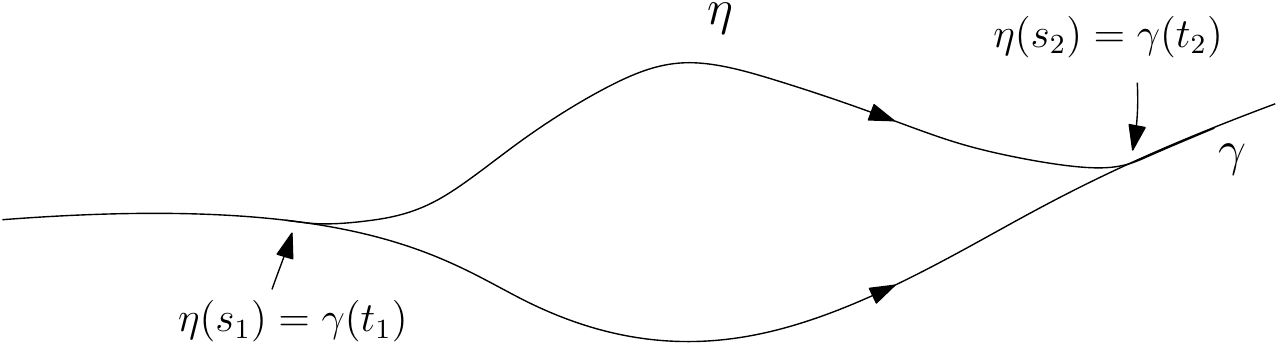}
\captionof{figure}{}
\end{center}

\medskip

$\gamma: \Rset\to \Rset^n$ is called {\it a universal global characteristics} if it is a characteristics for every viscosity solution to (\ref{cell}). 

For $p\in \Rset^n$, the collection of all universal characteristics
$$
\begin{array}{ll}
\widetilde {\mathcal{A}}_p&=\cup \{(\dot \gamma(t), \gamma(t))|\ t\in \Rset, \ \text{$\gamma$ is a universal global characteristics}\}\\[3mm]
&=\cap_{\text{$v$ is a solution to (\ref{cell})} }J_v
\end{array}
$$
is defined as the Aubry set. Hence the following graph property holds:
\be{}\label{graphsupport}
\widetilde {\mathcal{A}}_p \subset \{(q,x)\in \Rset ^n\times \Rset^n \,:\,
\text{$Dv(x)$ exists and $p+Dv(x)=D_{q}L(q,x)$}\}. 
\ee
Write $\mathcal {A}_p$ as the projection of $\tilde {\mathcal{A}}_p$ on $\Rset^n$. 

Also, we define the collection of all global characteristics associated with viscosity solutions of (\ref{cell})
$$
\begin{array}{ll}
\tilde {\mathcal{N}}_p&=\cup \{(\dot \gamma(t), \gamma(t))|\ t\in \Rset, \ \text{$\gamma$ is a global characteristics}\}\\[3mm]
&=\cup_{\text{$v$ is a solution to (\ref{cell})} }J_v
\end{array}
$$
as the Man\'e set. 

In standard definitions (\cite{F}), Aubry set and Man\'e set are on $\Rset^n\times \Bbb T^n$. Here for convenience, we lift them to $\Rset ^n\times \Rset^n$. Let us highlight several key properties

\medskip

\nit{\bf $\bullet$ Property 4 } Two absolutely minimizing curves can not intersect twice unless they are the same after suitable translations in time. This property together with 2d topology plays a crucial role in the classical Aubry-Mather theory in two dimensions. 

\medskip

\nit{\bf $\bullet$ Property 5} Two universal global characteristics can not intersect unless they are the same after suitable translation in $t$;

\medskip

\nit{\bf $\bullet$ Property 6 } Any global characteristics, when it is projected to $\Bbb T^n$, can not intersect itself unless the orbit is periodic.

\medskip

\nit{\bf $\bullet$ Property 7} If $\xi$ is a global characteristics, then for any sequence $T_m\to \infty$ as $m\to +\infty$, 
\be\label{rotation-vector}
\lim_{m\to +\infty}{\xi (T_m)-\xi(0)\over T_m}\in \partial \overline H(p)
\ee
if the limit exists. The full limit 
$$
\lim_{T\to \infty}{\xi (T)-\xi(0)\over T},
$$
if exists, is called the {\it rotation vector} of $\xi$. 

\subsection{Mather set}

Denote by $\mathcal{W} $ the set of all Borel probability
measures on $\Rset^n \times \Bbb T^n$ which are Euler-Lagrangian
flow invariant. For fixed $p\in \Rset^n$, $\mu\in \mathcal{W}$ is called a {\it ``Mather measure"} if 
$$
\int_{\Rset^n\times \Bbb T^n} (L(q,x)-p\cdot q)\,d\mu=\min_{\nu\in \mathcal{W}}\int_{\Rset^n\times \Bbb T^n}(L(q,x)-p\cdot q)\,d\nu.
$$
Denote by $\mathcal{W}_p$ the set of all such Mather measures. The value of the minimum action on the right hand side turns out to be $-\overline H(p)$, i.e., 
\be\label{negative-hbar}
\min_{\nu\in \mathcal{W}}\int_{\Rset^n\times \Bbb T^n}(L(q,x)-p\cdot q)\,d\nu=-\overline H(p).
\ee
In dynamical system literature, the effective Hamiltonian $\overline H$ is called ``{\it $\alpha$-function}" and is often denoted as $\alpha(c)$, where $c$ is the same as $p$.

The {\it Mather set} is defined to be the closure of the union of the support of all Mather measures, i.e., 
$$
\widetilde{\mathcal{M}}_p=\overline {\bigcup_{\mu \in \mathcal{W}_p}\mathrm{supp}(\mu)}.
$$
The projected Mather set $\mathcal {M}_{p}$ is the projection of $\widetilde{\mathcal {M}}_p$ to the torus. A curve $\xi: \Rset\to \Bbb T^n$ is called an orbit on $\mathcal {M}_p$ if it satisfies (\ref{E-L}) and 
$$
(\dot \xi(0), \xi(0))\in \widetilde{\mathcal {M}}_p.
$$ 
If we lift $\widetilde{\mathcal{M}}_p$ to $\Rset^n\times \Rset^n$ (or project $\widetilde{\mathcal{A}}_p$ and $\widetilde{\mathcal{N}}_p$ in our definition to $\Rset^n\times \Bbb T^n$) , the following relation holds
$$
\widetilde{\mathcal{M}}_p\subset \widetilde{\mathcal{A}}_p\subset\widetilde{\mathcal{N}}_p.
$$
In particular, the graph property (\ref{graphsupport}) also holds for $\tilde {\mathcal {M}}_p$. Also, all trajectories on $\mathcal {M}_p$ are universal global characteristics and, hence, absolutely minimizing curves with respect to $\overline H(p)$. Moreover, all viscosity solutions to the cell problem (\ref{cell}) are $C^{1,1}$ on $\mathcal {M}_p$. See \cite{EG} for instance.

One hope is that Mather sets might have some sort of ``integrable structure" in the sense that long term behaviors of trajectories there can be better understood. However, when $n\geq 3$, very little has been known in this direction except in certain special cases like the classical Hedlund example and its generalizations (\cite{H}, \cite{JTY}). We also would like mention that for generic V, Corollary 1.2 in \cite{BC} says that there are at most $n+1$ ergodic Mather measures for every $p$, which is proved under some framework of convex analysis.

\medskip

The following Lemma is a well known fact in the theory of Hamilton-Jacobi equations \cite{Lions}.

\begin{lem}\label{calib} Let $U$ be an open subset of $\Rset^n$. Assume that for some $c\in \Rset$, $w\in W^{1,\infty}(U)$ satisfies that 
$$
H(Dw,x)\leq c \quad \text{for a.e $x\in U$}.
$$
Then for any $\eta\in AC([t_1,t_2], U)$,
$$
\int_{t_1}^{t_2} L(\dot\eta(t),\eta (t))+c\,dt\geq w(\eta(t_2))-w(\eta(t_1)).
$$
The equality holds if and only if $\eta$ is a characteristics of $w$, i.e., $w$ is differentiable along $\eta$, $H(Dw(\eta(t)),\eta(t))=c$ and $Dw(\eta(t))=D_qL(\dot \eta(t), \eta(t))$ for $t\in [t_1,t_2]$. 
\end{lem}

Using solutions to the cell problem (\ref{cell}), as an immediate corollary, we have that 
\begin{cor}\label{p-calib} If $\gamma:[0,T]\to \Rset$ satisfies that 
$$
\gamma(T)-\gamma(0)=\vec{l}\in \Zset^n,
$$
then 
$$
\int_{0}^{T}L(\dot \gamma,\gamma)+\overline H(p)\,dt\geq p\cdot \vec{l}.
$$
The equality holds if and only if $\gamma$ is a periodic orbit in $\mathcal{M}_p$. 
\end{cor}

\subsection{The classical Aubry-Mather theory when $n=2$}

In this section, we assume that $n=2$ and focus on
$$
H(p,x)={1\over 2}|p|^2+V.
$$
Throughout this section, we fix 
$$
c>\max_{\Rset ^n}V.
$$

Due to (\ref{rotation-vector}), the 2d topology and fact that different trajectories on $\mathcal {M}_p$ can not intersect, it can be proved that the level curve (\cite{C})
$$
S_c=\{\overline H(p)=c\}
$$
is $C^1$. This is equivalent to saying that for every $p\in S_c$, there exists a unit vector $q_p$ such that 
\be\label{normal-vector}
\partial \overline H(p)=\{\lambda q_p|\ \ \lambda\in [a,b]\}.
\ee
Here $a\leq b$ are two positive constants depending on $p$. $q_p$ is the outward unit normal vector of $S_c$ at $p$. Note that if $a=b$, then $\overline H$ is differentiable at $p$.

\medskip

\nit{\bf $\bullet$ Property 8: Cornerstone of the Aubry-Mather theory.} If $q_p$ is a rational vector, then all orbits on $\mathcal {M}_p$ are periodic orbits with the same first homology class $(m,n)\in \Zset^2$ that is irreducible.

\begin{rmk}\label{extreme-orbit} Choose $p_k\to p$ as $k\to +\infty$ such that $\overline H(p_k)>\overline H(p)$ and $q_{p_k}=q_p$, then by standard convex analysis, 
$$
\lim_{k\to +\infty}\partial \overline H(p_k)=bq_p.
$$
Together with the stablity of periodic orbits (see the analysis in the appendix), we can deduce that there is a periodic orbit on $\mathcal {M}_p$ whose rotation vector is $bq_p$. Similarly, there is a periodic orbit on $\mathcal {M}_p$ whose rotation vector is $aq_p$. 
\end{rmk}

\nit{\bf $\bullet$ Property 9: Identification with circle homeomorphism.} Choose $\hat p\in S_c$ such that $q_{\hat p}=(0,1)$. Let $\xi$ be a periodic orbit on $\mathcal{M}_{\hat p}$ and lift it to $\Rset^2$. Now for each $k\in \Zset$, denote
$$
\xi_k=\xi+(k,0).
$$
For $p\in \Rset^2$ with $q_p\not= (0,1)$ or $(0,-1)$, let $\gamma:\Rset\to \Rset^2$ be a global characteristics associated with a solution to the cell problem (\ref{cell}). Owing to Property 1 in section 2, for each $k\in \Zset$, $\gamma$ intersects with $\xi_k$ exactly once. Let $a_k\in \Rset$ be such that 
\[
\gamma\cap \xi_k=\xi_k(a_kT).
\]
\begin{center}
\includegraphics[scale=0.4]{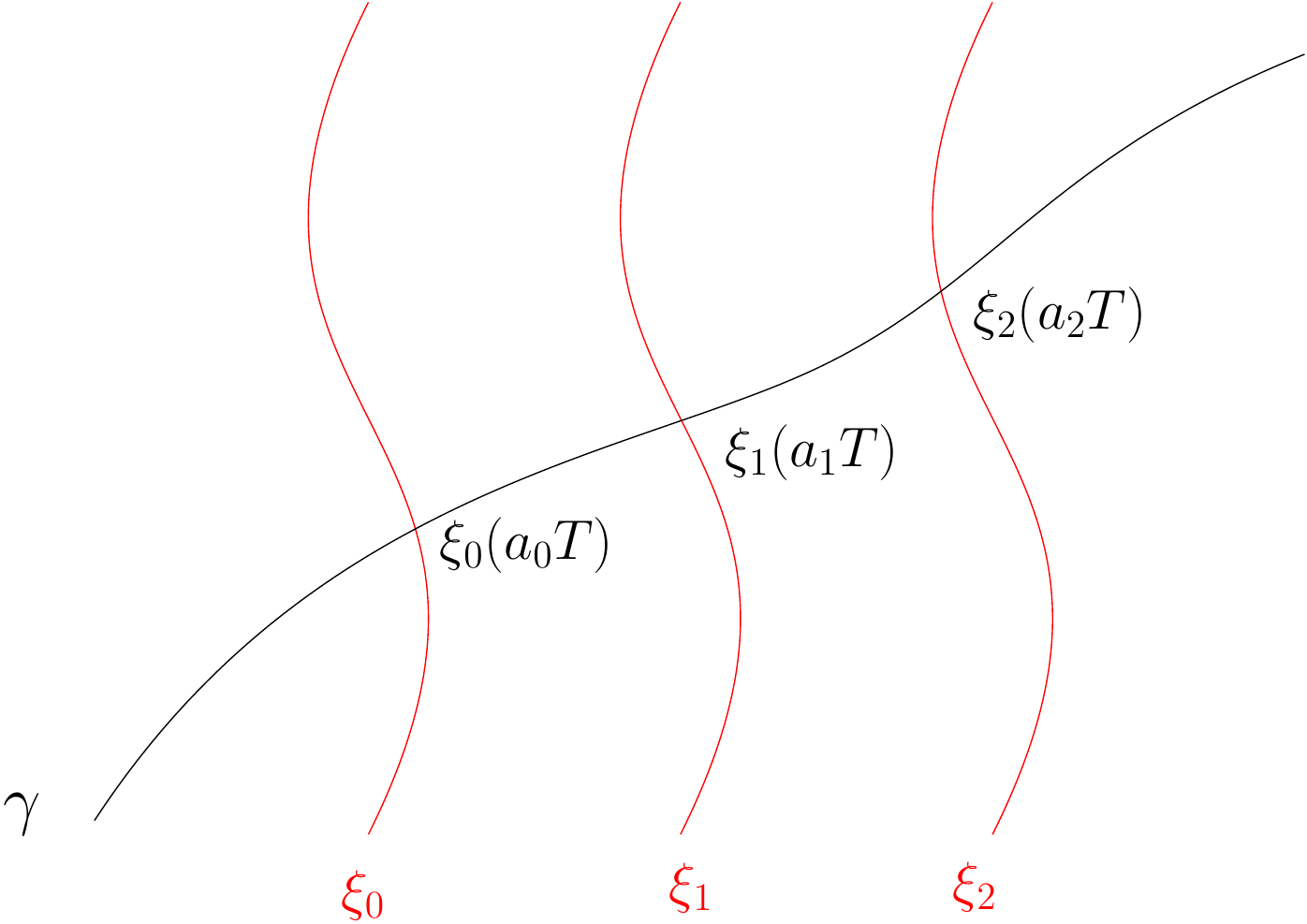}
\captionof{figure}{}
\end{center}

Since $\gamma$ can not intersect itself when it is projected to $\Bbb T^2$,
either $a_k=0$ for all $k\in \Zset$ or $\{a_k\}_{k\in \Zset}$ is a strictly monotonic sequence. If $\{a_k\}_{k\in \Zset}$ is strictly increasing, there exists a circle homeomorphism $f$ such that 
\[
f(a_k)=a_{k+1} \quad \text{for all $k\in \Zset$}.
\]
See \cite[Theorem 3.15]{B1} for further details on the definition of $f$. If $\{a_k\}_{k\in \Zset}$ is strictly decreasing, then consider $f(a_k)=a_{k-1}$. 

\medskip

We say that $p\in S_c$ is a {\it linear point} if there exists $p'\in \Rset^2$ such that $p'\not=p$ and the line segment (edge)
$$
[p,p']\subset S_c.
$$
Clearly, $q_p=q_{p'}$ and $q_p\cdot (p-p')=0$. Combining with the definition of Mather sets, we can deduce that, for every $p''\in [p,p']$,
\be\label{Matherequal}
\tilde {\mathcal{M}}_{p''}=\tilde{\mathcal{M}}_p.
\ee

The following two results were poved in \cite{B2}.

\medskip

\nit{\bf $\bullet$ Property 10} $p\in S_c$ is a linear point if and only if $\mathcal {M}_p\not=\Bbb T^2$. Moreover, $p$ is a linear point, then $q_p$ is a rational vector. The converse might not be true in general. Nevertheless, for generic V, Corollary 1.2 in \cite{BC} and the above Property 8 imply that if $q_p$ is an rational vector, then $\mathcal {M}_p$ has at most 3 periodic orbits and, hence, $p$ must be a linear point. See \cite{CZ} for hyperbolicity of periodic orbits. 

\medskip

\nit{\bf $\bullet$ Property 11} $S_c$ is not strictly convex (i.e., the set of linear points on $S_c$ is not empty) unless $V$ is a constant.

\medskip

Results in \cite{B2} are presented under the framework of minimizing geodesic associated with a periodic Riemannian metric on $\Rset^2$. Note that for mechanical Hamiltonians, every absolutely minimizing curve associated with $L(q,x)+c$ is a minimizing geodesics associated with the Riemannian metric
$$
g=\sqrt{2(c-V(x))}(dx_1\otimes dx_1+dx_2\otimes dx_2).
$$
The converse is also true after proper reparametrization.

\section{Proof of Theorem \ref{main2}}

Fix $p_0\in \Rset^2$ with $\overline H(p_0)>\max_{\Bbb T^2}V$. Denote 
$$
c_0=\overline H(p_0)
$$
and 
$$
S_{c_0}=\{p\in \Rset^2|\ \overline H(p)=c_0\}.
$$
Assume that the outward unit normal vector $q_{p_0}$ (see (\ref{normal-vector})) at $p_0$ is a rational vector. Denote
$$
q_{p_0}={(m,n)\over \sqrt{m^2+n^2}}.
$$
Here
$$
\text{$(m,n):$ the first homology class of periodic orbits on $\mathcal{M}_{p_0}$. }
$$
Throughout this section, we lift $\mathcal{M}_{p_0}$ to $\Rset^2$ and still denote it as $\mathcal{M}_{p_0}$. In addition, for convenience, the lift of a periodic orbit on $\mathcal{M}_{p_0}$ is also called a periodic orbit.

Suppose that $v$ is a viscosity solution to 
$$
{1\over 2}|p_0+Dv|^2+V=c_0.
$$
For $A, B\subset \Rset^2$ and 
$$
u=p_0\cdot x+v,
$$
we define the barrier between two sets as 
$$
d_{u}(A,B)=\inf_{x\in A, y\in B} (h(x,y)-(u(y)-u(x)))
$$
Here
$$
h(x,y)=\inf_{\substack{t>0,\ \gamma\in AC([0,t])\\ \gamma(0)=x,\ \gamma(t)=y}}\left( \int_{0}^{t}{1\over 2}|\dot \gamma|^2-V(\gamma)+c_0\,ds\right)
$$
and $AC([0,t])$ the set of all absolutely continuous curve $[0,t]\to \Rset^2$. Owing to Lemma \ref{calib}, 
$$
d_u(A,B)\geq 0.
$$

If $L_1:\Rset\to \Rset^2$ and $L_2:\Rset\to \Rset^2$ ar two curves, $d_u(L_1,L_2)$ is understood as $d_u(A_1, A_2)$ for $A_1=L_1(\Rset)$ and $A_2=L_2(\Rset)$.

\begin{lem} Suppose that $\gamma_1$ and $\gamma_2$ are two periodic orbits on $\mathcal{M}_{p_0}$. Then 
$$
d_u(\gamma_1,\gamma_2)=\lim_{t\to +\infty, \ s\to -\infty }(h(\gamma_1(s),\gamma_2(t))-(u(\gamma_2(t))-u(\gamma_1(s)))).
$$
\end{lem}

\nit Proof: Denote 
$$
G(t,s)=h(\gamma_1(s),\gamma_2(t))-(u(\gamma_2(t))-u(\gamma_1(s))).
$$
We claim that $G$ is decreasing on $t$ and increasing $s$. In fact, for $t_1<t_2$, due to the triangle inequality, 
$$
\begin{array}{ll}
h(\gamma_1(s),\gamma_2(t_2))&\leq h(\gamma_1(s),\gamma_2(t_1))+h(\gamma_2(t_1),\gamma_2(t_2))\\[5mm]
&= h(\gamma_1(s),\gamma_2(t_1))+(u(\gamma_2(t_2))-u(\gamma_2(t_1))))
\end{array}
$$
Hence $G(s,t_2)\leq G(s,t_1)$. Similarly, we can show that $G$ is increasing on $s$. Then our lemma follows immediately. \qed

Let $D$ be an open set bounded by two unbounded simple curves $L_1$ and $L_2$ on $\Rset ^2$. For any $\delta$ with $0\leq \delta \leq d_u(L_1, L_2)$, let
$$
\begin{cases}
g(x)=u(x) \quad \text{for $x\in L_1$}\\
g(x)=u(x)+\delta \quad \text{for $x\in L_2$}
\end{cases}
$$
\begin{center}
\includegraphics[scale=0.5]{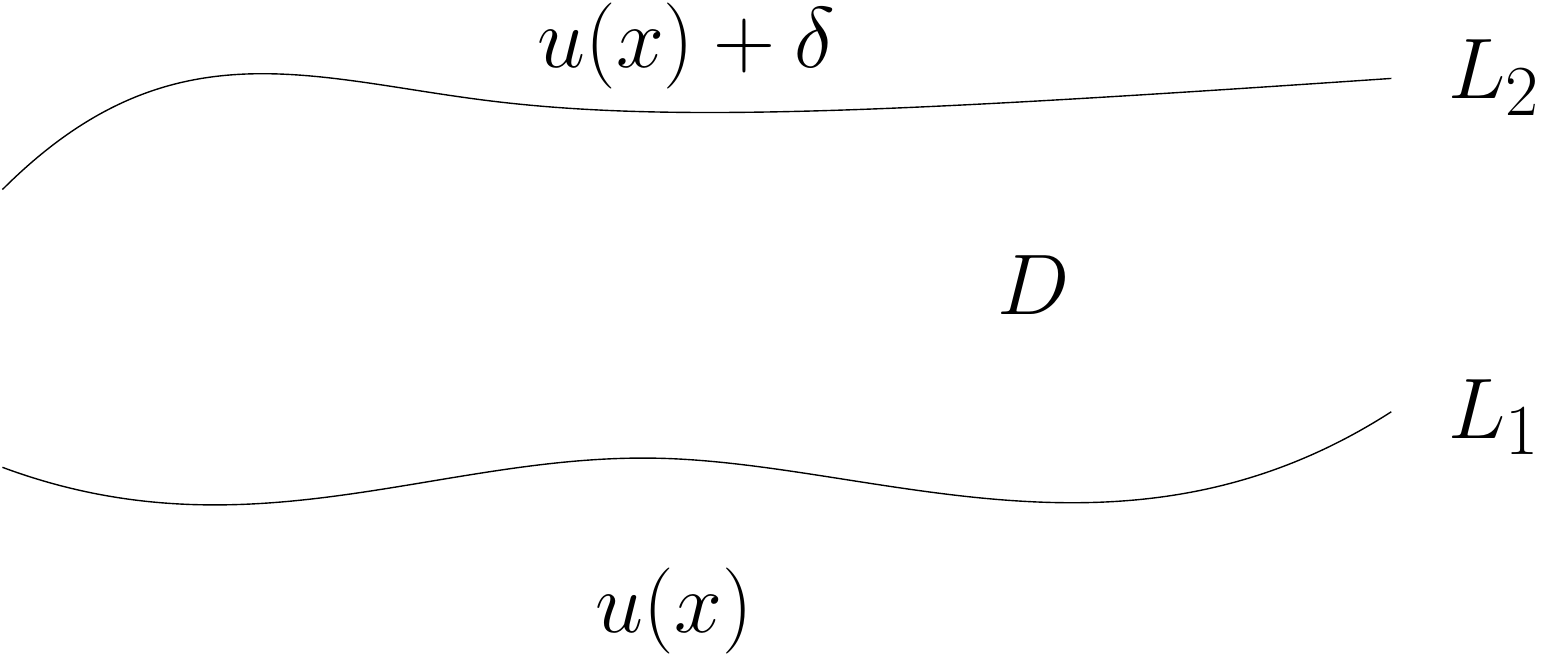}
\captionof{figure}{}
\end{center}
Define
$$
u_\delta(x)=\inf_{y\in L_1\cup L_2}\{g(y)+h(y,x)\} \quad \text{for $x\in \bar D$}.
$$
The following lemma is due to the well known compatibility condition for existence of solutions to Hamilton-Jacobi equations \cite{Lions}.
\begin{lem}\label{compatibility}
$u_\delta$ is a Lipschitz continuous viscosity solution to
$$
\begin{cases}
{1\over 2}|Du_{\delta}|^2+V(x)=c_0\\
u_{\delta}=u \quad \text{on $L_1$}\\
u_{\delta}=u+\delta \quad \text{on $L_2$}.
\end{cases}
$$
\end{lem}
As a corollary, we have that 
\begin{cor}\label{no-gap} Let $\xi_1$ and $\xi_2$ be two periodic orbits on $\mathcal {M}_{p_0}$. If the region $D\subset \Rset^2$ bounded by $\xi_1$ and $\xi_2$ is foliated by periodic orbits on $\mathcal {M}_{p_0}$, then
$$
d_u(\xi_1, \xi_2)=0.
$$
\begin{center}
\includegraphics[scale=0.3]{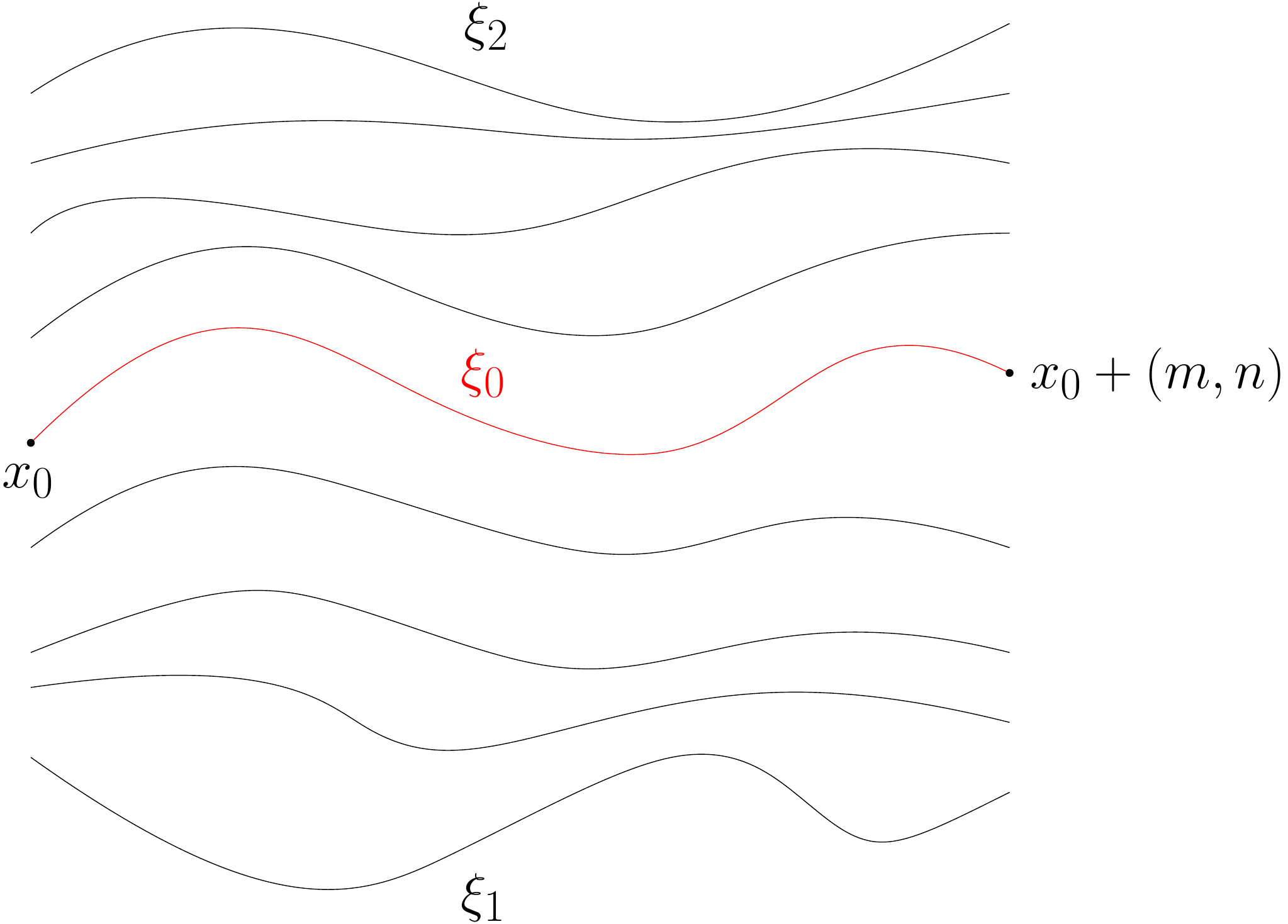}
\captionof{figure}{}
\end{center}
\end{cor}

Proof: Let $(m.n)$ be the homology class of periodic orbits on $\mathcal {M}_p$. Choose $\delta=d_u(\xi_1, \xi_2)$. Apparently, $D+(m,n)=D$. Also,
$$
\begin{array}{ll}
u_\delta (x+(m,n))&=\inf_{y\in L_1\cup L_2}\{g(y)+h(y,x+(m,n))\}\\[3mm]
&=\inf_{y\in L_1\cup L_2}\{g(y+(m,n))+h(y+(m,n),x+(m,n))\}\\[3mm]
&=\inf_{y\in L_1\cup L_2}\{g(y+(m,n))+h(y,x)\}\\[3mm]
&=\inf_{y\in L_1\cup L_2}\{g(y)+h(y,x)\}+p_0\cdot (m,n)\\[3mm]
&=u_\delta(x)+p_0\cdot (m,n).
\end{array}
$$
The second equality in the above is due to $(L_1\cup L_2)+(m,n)=L_1\cup L_2$. The third one is owing to the periodicity of $V$. The last equality is because $g(y+h)-g(y)=u(y+h)-u(y)=p_0\cdot h$ for all $h\in \Rset^2$. 

To prove $d_u(\xi_1, \xi_2)=0$, it suffices to show that 
$$
u_{\delta}\equiv u \quad \text{in $\bar D$}.
$$
In fact, choose an arbitrary point $x_0\in D$ and let $\xi_0$ be the periodic orbit that passes through $x_0$ with minimum period $T_0$. Then 
$$
\begin{array}{ll}
u_{\delta}(x_0+(m,n))-u_{\delta}(x_0)=(m,n)\cdot p_0&=u(x_0+(m,n))-u(x_0)\\[5mm]
&=\int_{0}^{T_0}{1\over 2}|\dot \xi_0|^2-V(\xi_0)+c_0\,ds.
\end{array}
$$
The second ``=" is due to the fact that every trajectory on $\mathcal{M}_{p_0}$ is a universal global characteristics. Hence $\xi_0$ is also a characteristics of $u_{\delta}$ and $Du_{\delta}(x_0)=Du(x_0)=\dot \xi_0(0)$. Accordingly, 
$$
Du_{\delta}\equiv Du \quad \text{in $D$}.
$$
Therefore $u_{\delta}=u$ in $\bar D$. \qed

\begin{lem}\label{flat-condition} Suppose that $\xi$ is a periodic orbit on $\mathcal{M}_{p_0}$. Let $\xi_2=\xi_1+(-n,m)$. 

(1) There exists $\tau>0$ such that 
$$
E_{p_0,\tau,+}=\{p_0+t(-n,m)|\ t\in [0, \tau]\}\subset S_{c_0},
$$
if and only if
$$
d_u(\xi_1, \xi_2)>0.
$$

(2) There exists $\tau>0$ such that 
$$
E_{p_0,\tau,-}=\{p_0+t(n,-m)|\ t\in [0, \tau]\}\subset S_{c_0},
$$
if and only if 
$$
d_u(\xi_2, \xi_1)>0.
$$
\end{lem}

Proof: It is enough to prove (1). The proof for (2) is similar. 

\medskip

$``\Rightarrow"$ Denote 
$$
p_{\tau}=p_0+\tau (-n,m).
$$
Let $v_\tau$ be a viscosity solution to 
$$
{1\over 2}|p_{\tau}+Dv_{\tau}|^2+V(x)=c_0.
$$
Recall that $c_0=\overline H(p_0)$. Owing to (\ref{Matherequal}), $\xi_1$ and $\xi_2$ are also periodic orbits on $\mathcal {M}_{p_\tau}$. Hence 
$$
Du_{\tau}=Du \quad \text{on $\xi_1\cup \xi_2$}
$$
Without loss of generality, we may assume that $u_\tau=u$ on $\xi_1$. Then $u_\tau=u+\tau(m^2+n^2)$ on $\xi_2=\xi_1+(-n,m)$. Accordingly,
$$
d_u(\xi_1,\xi_2)=d_{u_\tau}(\xi_1,\xi_2)+\tau (m^2+n^2)\geq \tau (m^2+n^2)>0.
$$

\medskip

$``\Leftarrow"$ The part is not really needed for our main result. Since it is the essential PDE part in proving the existence of edges of $S_{c_0}$, we present it here. Choose 
$$
\tau={d_u(\xi_1,\xi_2)\over m^2+n^2}.
$$ 
For $p'\in E_{p_0,\tau,+}$, define
$$
F=
\begin{cases}
u \quad \text{on $\xi_1$}\\
u+p'\cdot (-n,m) \quad \text{on $\xi_2$}
\end{cases}
$$
Clearly, 
$$
\inf_{x,y\in \xi_1\cup \xi_2} (h(x,y)-(F(y)-F(x)))\geq 0.
$$
Hence for $x\in D$ 
$$
u_{p'}(x)=\inf_{y\in \xi_1\cup \xi_2} \{F(y)+h(x,y)\}
$$
is a viscosity solution of 
$$
\begin{cases}
{1\over 2}|Du_{p'}|^2+V(x)=c_0 \quad \text{in $D$}\\
u_{p'}=F \quad \text{on $\partial D$}.
\end{cases}
$$
As in Corollary \ref{no-gap}, we have that for $x\in D$,
$$
u_{p'}(x+(m,n))=u_{p'}(x)+p_0\cdot (m,n)=u_{p'}(x)+p'\cdot (m,n).
$$
Next we may extend $u_{p'}$ to $\Rset^2$ by
$$
u_{p'}(x+k(-n,m))=u_{p'}(x)+kp'\cdot (-n,m) \quad \text{for all $x\in D$ and $k\in \Zset$}.
$$
See Figure 5 below. Apparently, $v_{p'}=u_{p'}-p'\cdot x$ is $(m^2+n^2)\Zset^2$ periodic viscosity solution to
$$
{1\over 2}|p'+Dv_{p'}|^2+V=c_0 \quad \text{in $\Rset^2/(\xi_1+(-n,m)\Zset)$}.
$$
Since $\xi_1$ is a periodic characteristics of $v_{p'}$, it is easy to see that $v_{p'}$ is a $(m^2+n^2)\Zset^2$ periodic viscosity solution to
$$
{1\over 2}|p'+Dv_{p'}|^2+V=c_0 \quad \text{on $\Rset ^2$}.
$$
By the inf-max formula (\ref{inf-max}), it is not hard to deduce that 
$$
\overline H(p')=c_0.
$$
\begin{center}
\includegraphics[scale=0.3]{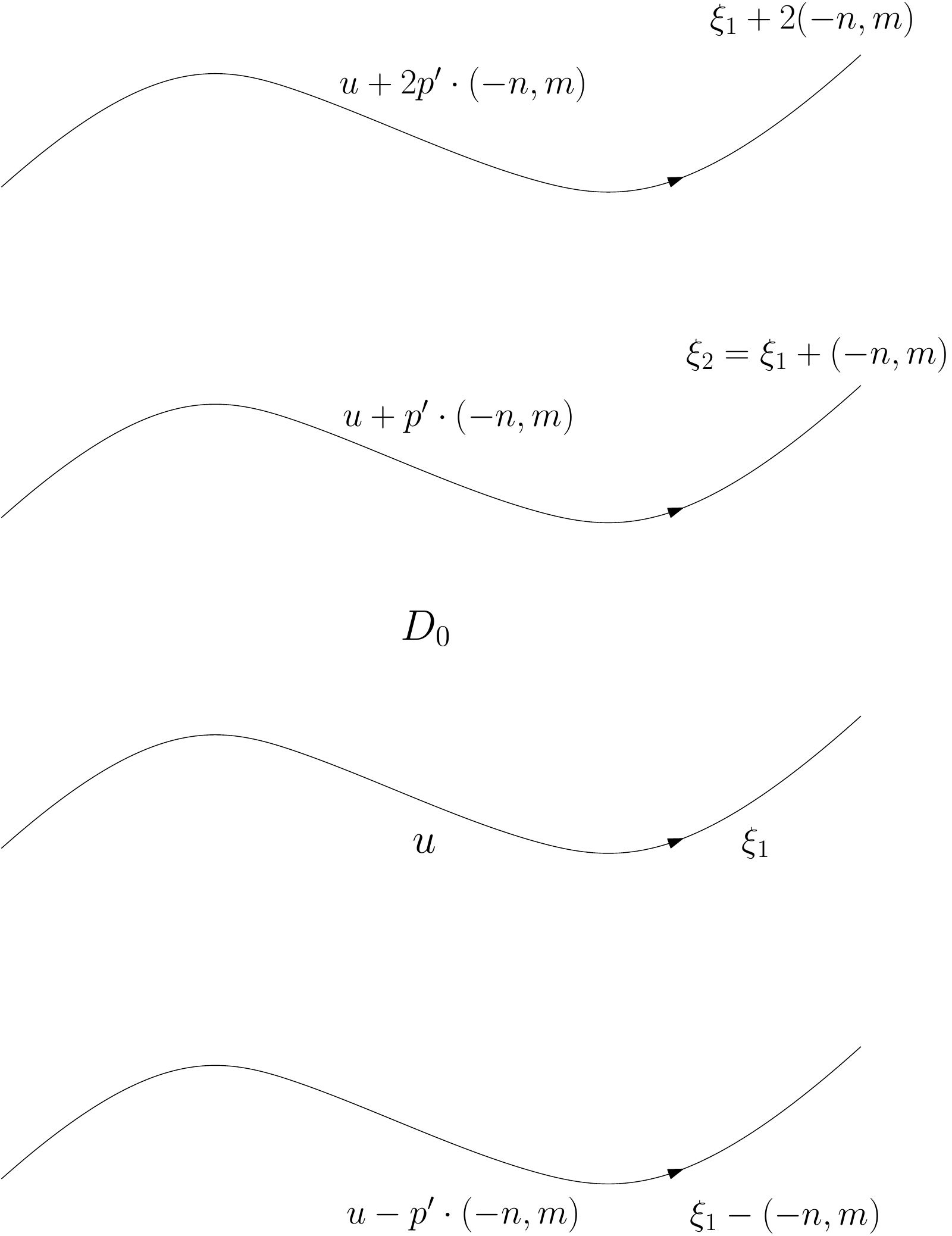}
\captionof{figure}{}
\end{center}
\qed
\begin{lem}\label{curve-linear} Suppose that $L_1$, $L_2$, $L_3$, $\cdots$, $L_m$ are $m$ different periodic orbits on $\mathcal {M}_{p_0}$. For $i=2,\cdots m-1$, $L_i$ lies between $L_{i-1}$ and $L_{i+1}$. Then
$$
d_u(L_1,L_m)=\sum_{k=1}^{m-1}d_u(L_k,L_{k+1}).
$$
\end{lem}

Proof: By induction, it suffices to establish the above equality for $m=3$. By definition of $d_u$, it is clear that 
$$
d_u(L_1, L_3)\geq d_u(L_1,L_2)+d_u(L_2, L_3).
$$
Now let us prove the other direction. According to the definition, for any $\delta>0$, there exist two Lipschitiz continuous curves $\gamma_1: [0,a]\to \Rset ^2$ and $\gamma_2: [0,b]\to \Rset^2$ such that
$$
\gamma_1(0)\in L_1, \quad \gamma_1(a),\ \gamma_2(0)\in L_2, \quad \gamma_2(b)\in L_3
$$
$$
d_u(L_1,L_2)\geq \int_{0}^{a}{1\over 2}|\dot \gamma_1|^2-V(\gamma_1)+c_0\,ds-u(\gamma_1(0))+u(\gamma_1(a))-\delta.
$$
and
$$
d_u(L_2,L_3)\geq \int_{0}^{b}{1\over 2}|\dot \gamma_2|^2-V(\gamma_2)+c_0\,ds-u(\gamma_2(0))+u(\gamma_2(b))-\delta.
$$
By periodic translation, we can make $\gamma_1(a)$ at an earlier time than $\gamma_2(0)$, i.e., there exists $t_1<t_2$ such that 
$$
\gamma_1(a)=L_2(t_1) \quad \mathrm{and} \quad \gamma_2(0)=L_2(t_2).
$$
By connecting $\gamma_1$, $L_2$ and $\gamma_2$, we define $\gamma: [0, \ a+t_2-t_1+b]\to \Rset^2$ as
$$
\gamma (t)=
\begin{cases}
\gamma_1(t) \quad \text{for $t\in [0,a]$}\\
L_2(t-a+t_1) \quad \text{for $t\in [a, a+t_2-t_1]$}\\
\gamma_2(t+t_t-t_2-a) \quad \text{for $t\in [a+t_2-t_1, a+t_2-t_1+b]$}.
\end{cases}
$$
\begin{center}
\includegraphics[scale=0.5]{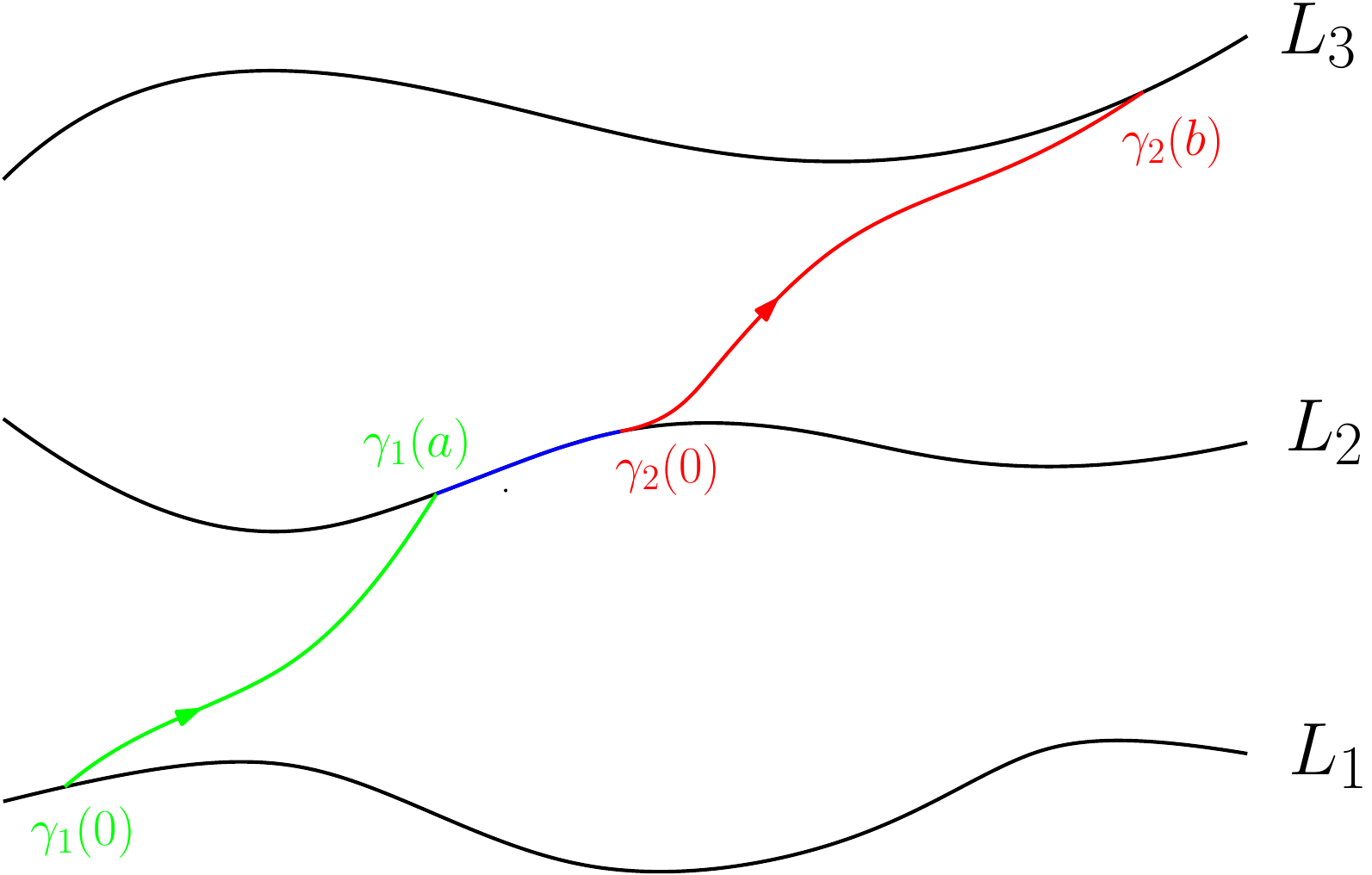}
\end{center}
Since
$$
\int_{t_1}^{t_2}{1\over 2}|\dot L_2 |^2-V(L_2)+c_0\,ds-u(L_2(t_2))+u(L_2(t_1))=0,
$$
we have that for $\bar t=a+b+t_2-t_1$,
$$
\int_{0}^{\bar t}{1\over 2}|\dot \gamma |^2-V(\gamma)+c_0\,ds-u(\gamma(\bar t))+u(\gamma(0))\leq d_u(L_1,L_2)+d_u(L_2,L_3)+2\delta.
$$
Accordingly,
$$
d_u(L_1,L_3)\leq d_u(L_1,L_2)+d_u(L_2,L_3)+2\delta.
$$
Sending $\delta\to 0$ leads to 
$$
d_u(L_1,L_3)\leq d_u(L_1,L_2)+d_u(L_2,L_3).
$$
Hence our lemma holds. \qed

Next we prove that $\overline H$ is 1d near any point that is in the interior of an edge. This conclusion holds for any $V\in C^{k}(\Bbb T^2)$ when $k\geq 2$.

\begin{lem}\label{main1} Suppose there exist $p_1\not= p_2\in S_{c_0}$ such that 
$$
p_0\in \{tp_1+(1-t)p_2|\ t\in (0,1)\}\subset S_{c_0},
$$
i.e., $p_0$ is in the interior of an edge. Then there exists $r>0$ and a convex function $f:\Rset\to \Rset^2$ such that 
$$
\overline H(p)=f(q_{p_0}\cdot p) \quad \text{for $p\in B_r(p_0)$}.
$$
\end{lem}

Proof: It suffices to show that there exists $r>0$ such that
\be\label{vectorsame}
q_{p}= q_{p_0} \quad \text{for $p\in B_r(p_0)$}.
\ee
Recall the definition of $q_p$ in (\ref{normal-vector}). If this holds, $\overline H$ is constant along the direction that is perpendicular to $q_{p_0}$ in $B_r(p_0)$. Then
$$
f(t)=\overline H(p_0+q_{p_0}(t-p_0\cdot q_{p_0})).
$$
We argue by contradiction. If not true, then there exist $\{p_k\}_{k\geq 1}$ such that $p_k\to p_0$ as $k\to +\infty$ and 
$$
q_{p_k}\cdot (-n,m)\not=0.
$$
Without loss of generality, let us assume that
\be\label{positive-direction}
q_{p_k}\cdot (-n,m)>0.
\ee
Let $v_k$ be a periodic viscosity solution to 
$$
{1\over 2}|p_k+Dv_k|^2+V(x)=\overline H(p_k).
$$
subject to $\int_{\Bbb T^2}v_k\,dx=0$. Up to a subsequence if necessary, we assume that 
$$
\lim_{k\to +\infty}v_k=\hat v \quad \text{uniformly on $\Rset^2$}.
$$
Then $\hat v$ is a $\Zset^2$-periodic viscosity solution of
$$
{1\over 2}|p_0+D\hat v|^2+V(x)=c_0 \quad \text{on $\Rset^2$}
$$
subject to $\int_{\Bbb T^2}\hat v\,dx=0$. 

For $k\geq 1$, choose $\tilde p_{k}$ such that 
$$
\overline H(p_k)=\overline H(\tilde p_k)
$$
and
$$
q_{\tilde p_k}={(-n,m)\over \sqrt{m^2+n^2}}.
$$
Let $\eta_k:\Rset\to \Rset^2$ be a periodic orbit on $\mathcal{M}_{\tilde p_k}$ (lift to $\Rset^2$). Up to a subsequence, we may assume that
$$
\quad \lim_{k\to +\infty}\tilde p_k=\tilde p,
$$
and
$$
\lim_{k\to +\infty}\eta_k=\eta \quad \text{locally uniformly on $\Rset$}.
$$
Then $\overline H(\tilde p)=\overline H(p_0)$ and $\eta$ is a periodic orbit on $\mathcal {M}_{\tilde p}$ with first homology class $(-n,m)$. 
\begin{center}
\includegraphics[scale=0.5]{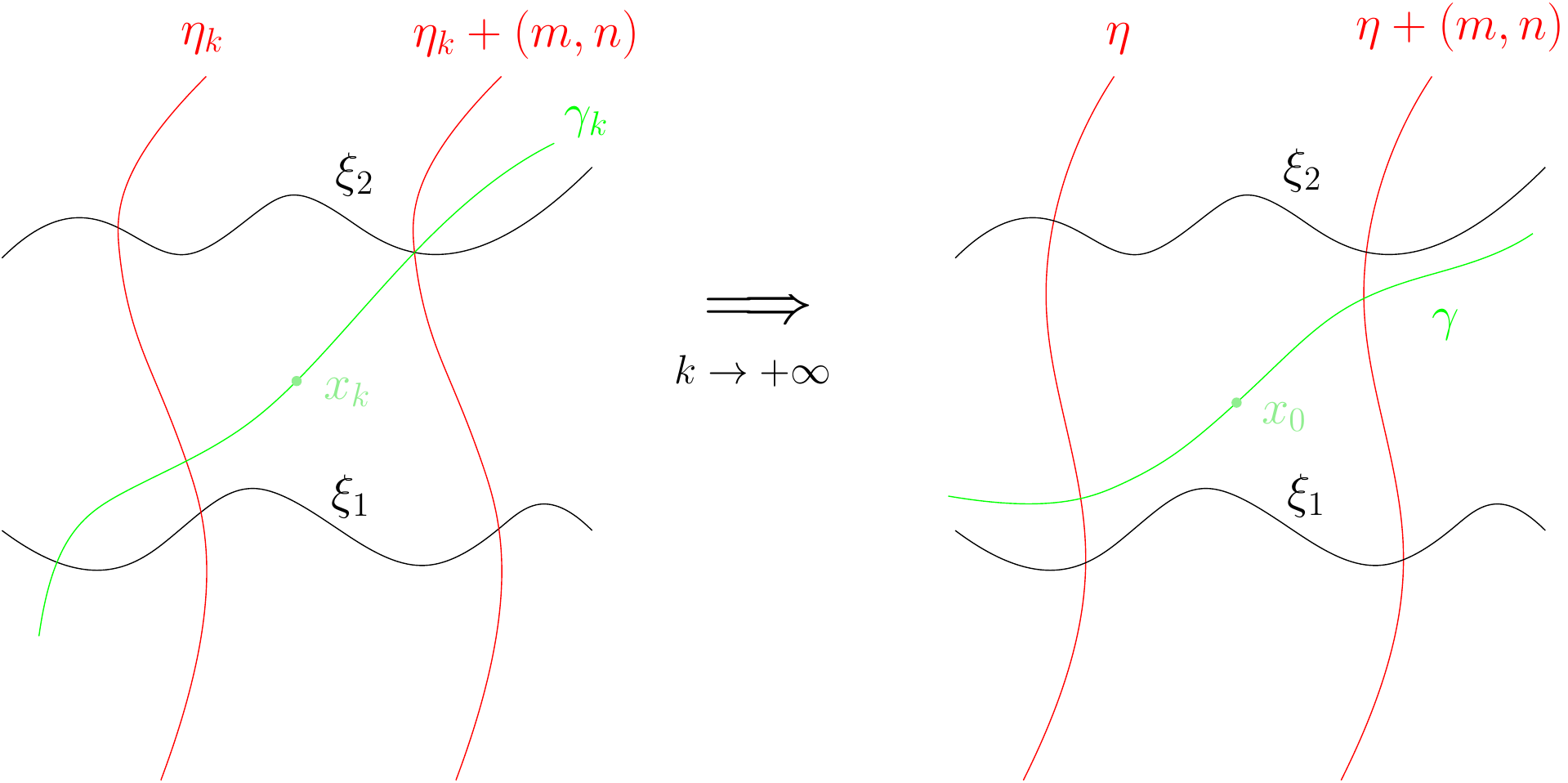}
\end{center}
\medskip

Suppose that $\xi_1$ and $\xi_2$ are two periodic orbits on $\mathcal{M}_{p_0}$ such that there is no other periodic orbit in the region bounded by $\xi_1$ and $\xi_2$ (a gap). 

\medskip

\nit{\bf Claim:} For $\hat u=p_0\cdot x+\hat v$
$$
d_{\hat u}(\xi_1,\xi_2)=0.
$$
$\xi_1$ divides the plane into two regions. Without loss of generality, we may assume that $\xi_2$ is on the same region as $\xi_1+(-n, m)$. For $k\geq 1$, let $\gamma_k:\Rset\to \Rset^2$ be an orbit on $\mathcal {M}_{p_k}$ (lift to $\Rset^2$) with $\gamma_k(0)=x_k$ for some $x_k\in [0,1]^2$ satisfying that the distance from $x_k$ to $\xi_2$ is half of the distance between $\xi_1$ and $\xi_2$:
$$
\mathrm{distance}(x_k,\xi_2)={1\over 2}\mathrm{distance}(\xi_1,\xi_2).
$$
Up to a subsequence if necessary, we assume that 
$$
\lim_{k\to +\infty}\gamma_k=\gamma \quad \text{locally uniformly on $\Rset$}.
$$
and
$$
\lim_{k\to +\infty}x_k=\bar x\in [0,1]^2.
$$
Due to the stability, $\gamma: \Rset\to \Rset^2$ is a global characteristics of $\hat v$ with $\gamma(0)=\bar x$. Owing to (\ref{positive-direction}), intersections of $\gamma$ with respect to $\eta+(m,n)\Zset$ are non-decreasing (see 6 in section 2.3 for the precise meaning). Since there is no periodic orbit between $\xi_1$ and $\xi_2$, there exist $\{t_{i}^{+}\}_{i\geq 1}$ and $\{t_{k}^{-}\}_{k\geq 1}$ such that 
$$
\lim_{i\to +\infty}t_{i}^{+}=+\infty \quad \mathrm{and} \quad \lim_{i\to +\infty}|\gamma (t_{i}^{+})-\xi_2 (t_{i}^{+})|=0
$$
and
$$
\lim_{i\to +\infty}t_{i}^{-}=-\infty \quad \mathrm{and} \lim_{i\to +\infty}|\gamma (t_{i}^{-})-\xi_1 (t_{i}^{-})|=0.
$$
Accordingly, we can derive that
$$
d_{\hat u}(\xi_1,\xi_2)=0
$$

Combining Corollary \ref{no-gap} and Lemma \ref{curve-linear}, we deduce that
$$
d_{\hat u}(\xi_1, \xi_1+(-n,m))=0.
$$
This contradicts to Lemma \ref{flat-condition} and the assumption that $p_0$ is in the interior of an edge. Hence (\ref{vectorsame}) and the Lemma holds. \qed

\medskip

\nit {\bf Proof of Theorem \ref{main2}.} Let 
$$
\mathcal {G}=\cap_{q\in \Qset^2, \ r\in \Qset}C(q,r).
$$
Here $\Qset^2$ is the collection of points on $\Rset^2$ whose both coordinates are rational numbers. See the appendix for the definition of $C(q,r)$. Fix $V\in \mathcal {G}$. Then for any $p\in \Rset^2$, if $\overline H(p)>\max_{\Rset^2}V$ and $\partial \overline H(p)\cap \Qset^2\not=\emptyset$, $p$ must be a linear point. Now write 
$$
O_V=\{p\in \Rset^2|\ \overline H(p)>\min_{\Rset^2}V \ \text{ and $p$ is an interior linear point}\}\cup F_{0}^{\circ}.
$$
Recall that the minimum level set $F_0=\{p\in \Rset^2|\ \overline H(p)=\max_{\Rset^2}V\}$. Owing to Lemma \ref{main1}, $O_V$ is an open set. We only need to show that $O_V$ is dense. Let us argue by contradiction. If not, then there exists $\bar p\in \Rset^2$ and $B_r(\bar p)$ for some $r>0$ such that $\overline H$ is strictly convex in $B_r(\bar p)$. Then the set
$$
W=\cup_{p\in B_r(\bar p)}\partial \overline H(p)
$$
is a non-empty open set. In particular, $W\cap \Qset^2$ is not empty. This contradicts the choice of $V$ and $O_V$. \qed

\section{Appendix}

\begin{defin}Given $r>0$ and a non-zero rational vector $q$. Denote by
$$
C(q,r)
$$
the collection of all $V\in C^{k}(\Bbb T^2)$ such that for any $p\in \Rset^2$, if 
$$
\overline H_V(p)\geq \max_{\Bbb T^2}V+r, \quad \mathrm{and} \quad q\in \partial \overline H(p)
$$
then
$$
\mathcal{M}_p\not=\Bbb T^2 \quad \text{(equivalently, $p$ is a linear point).}
$$
\end{defin}

\begin{theo} For $r>0$ and any non-zero rational vector $q$, $$
C(q,r)
$$
is an open dense set. 
\end{theo}

Proof: It is equivalent to showing that the complement 
$$
S=C^{k}(\Bbb T^2)\backslash C(q,r)
$$
is closed and nowhere dense. Write
$$
q=\lambda (m,n),
$$
where $\lambda>0$ and $(m,n)\in \Zset^2$ is irreducible. 

\medskip

{\bf Step 1:} We first prove that $S$ is closed. Suppose that $\{V_j\}_{j\geq 1}$ is a sequence of functions in $S$ and 
$$
\lim_{j\to +\infty}V_j=V \quad \text{in $C^{k}(\Bbb T^2)$}.
$$
Then
$$
\lim_{j\to +\infty}\overline H_{V_j}(p)=\overline H_V(p) \quad \text{locally uniformly in $\Rset^2$}.
$$
Since $V_j\in S$, for each $j\in \Nset$, there exists $p_j\in \Rset^2$ such that 
$$
\overline H_{V_j}(p_j)\geq \max_{\Bbb T^2}V+r, \qquad q\in \partial \overline H_{V_j}(p_j)
$$
and
$$
\mathcal{M}_{p_j,V_j}=\Bbb T^2.
$$
Due to the convexity, $ \overline H_{V_j}(0)\geq \overline H_{V_j}(p_j)+q\cdot (-p_j)$. Together with the quadratic growth of $\overline H_{V_j}$, it is easy to see that $\{p_j\}_{j\geq 1}$ is uniformly bounded. Upon a subsequence if necessary, we assume that 
$$
\lim_{j\to +\infty}p_j=\tilde p.
$$
Then 
$$
\overline H_V(\tilde p)\geq \max_{\Bbb T^2}V+r
$$
and by upper-semi-continuity of subdifferentials
$$
q\in \partial \overline H_V(\tilde p).
$$
Now we just need to show that
$$
\mathcal{M}_{\tilde p,V}=\Bbb T^2.
$$
This follows easily from stability of periodic orbits. For reader's convenience, we present details here. Fix $x_0\in \Bbb T^2$. 

For each $j\geq 1$, let $\xi_j: [0,T_j]\to \Bbb T^2$ be a periodic orbit on $\mathcal{M}_{p_j,V_j}$ with $\xi_j(0)=x_0$ and $T_j$ is the minimum period. Then
$$
\xi_j(T_j)=x_0+(m,n), \quad { \xi(T_j)-x_0\over T_j}\in \partial \overline H_{V_j}(p_j)
$$
and
$$
p_j\cdot (m,n)=\int_{0}^{T_j}{1\over 2}|\dot \xi_j|^2-V_j(\xi_j)+\overline H_{V_j}(p_j)\,dt.
$$
It is easy to see that $T_j$ and $||\xi_j||_{C^2([0,T_j])}$ are uniformly bounded. Up to a subsequence if necessary, we may assume that 
$$
\lim_{j\to +\infty}T_j=T \quad \mathrm{and} \quad \lim_{j\to +\infty}\xi_j=\xi \quad \text{uniformly in $C^1(\Rset)$}.
$$
Then
$$
\xi(0)=x_0 \quad \mathrm{and} \quad \xi(T)=x_0+(m,n)
$$
and
$$
\tilde p\cdot (m,n)=\int_{0}^{T}{1\over 2}|\dot \xi|^2-V(\xi)+\overline H_{V}(\tilde p)\,dt.
$$
So $\xi$ is a periodic orbit on $\mathcal {M}_{\tilde p,V}$ and $x_0\in \mathcal {M}_{\tilde p,V}$. Hence $\mathcal {M}_{\tilde p,V}=\Bbb T^2$. 

\medskip

{\bf Step 2:} Next we prove that $C(q,r)$ is dense. Given $V_0\in S$. Then there exists $p_0\in \Rset^2$ such that $\overline H_{V_0}(p_0)\geq \max_{\Bbb T^2}V_0+r$, $q\in \partial \overline {H}(p_0)$ and $\mathcal{M}_{p_0,V_0}=\Bbb T^2$. Assume that 
$$
\partial \overline {H}_{V_0}(p_0)=[\alpha q, \beta q]
$$
for $0<\alpha\leq 1\leq \beta$. By Remark \ref{extreme-orbit}, we can choose two periodic orbits $\xi_\alpha$ and $\xi_\beta$ on $\mathcal{M}_{p_0,V_0}$ such that their rotation vectors are $\alpha q$ and $\beta q$ respectively.

Choose $x_0\in \Bbb T^2$ such that $\xi_\alpha$ and $\xi_\beta$ do not pass through $x_0$. Pick $\delta>0$ such that 
$$
B_{2\delta}(x_0)\cap \xi_\alpha(\Rset)=\emptyset \quad \mathrm{and} \quad B_{2\delta}(x_0)\cap \xi_\beta(\Rset)=\emptyset.
$$
\begin{center}
\includegraphics[scale=0.5]{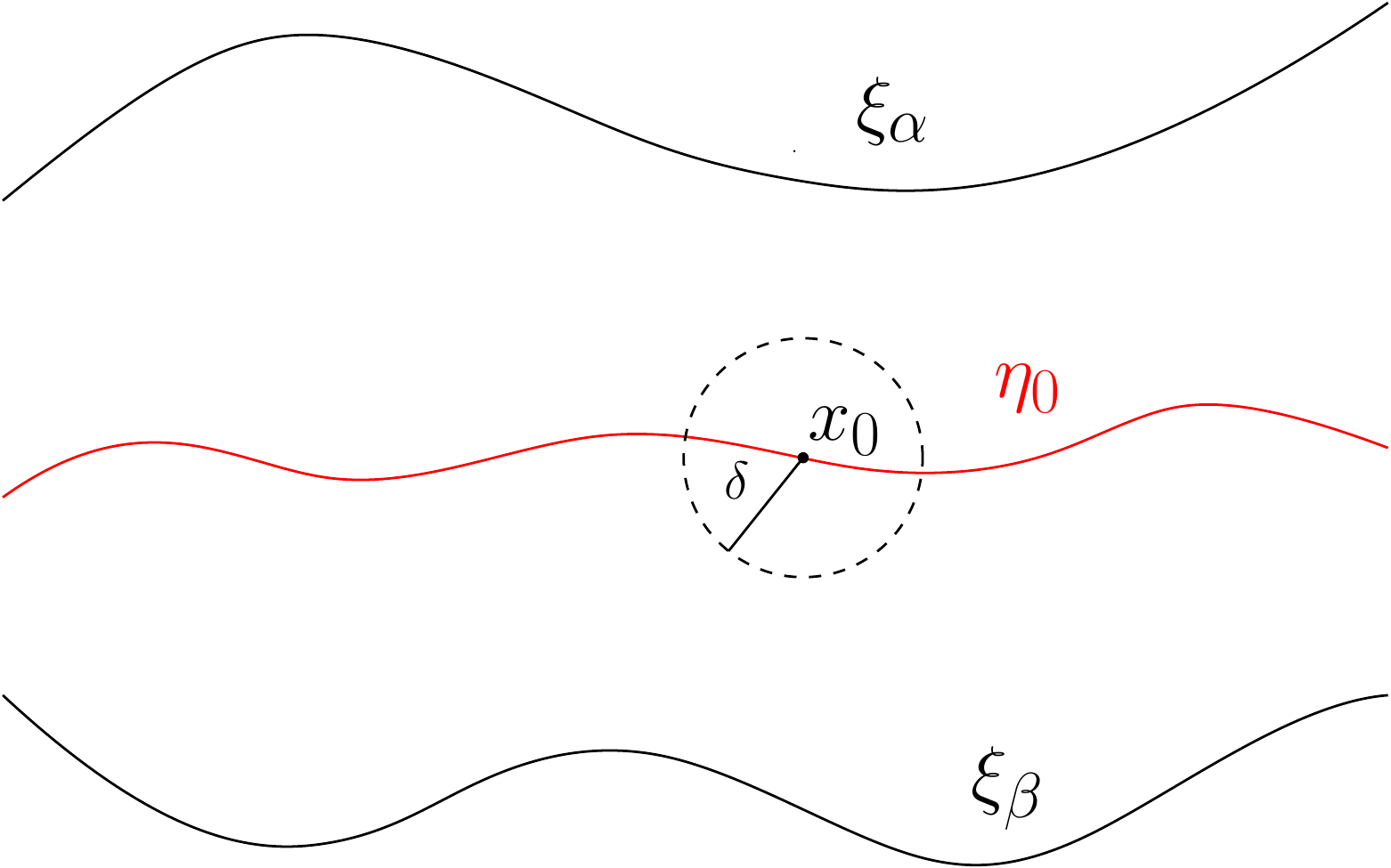}
\end{center}
Choose $\phi\in C^{\infty}(\Bbb T^2)$ satisfying that 
$$
\phi>0 \quad \text{in $B_\delta(x_0)$} \quad \mathrm{and} \quad \phi=0 \quad \text{in $\Bbb T^2\backslash B_\delta(x_0)$} 
$$
For $\delta>0$, denote
$$
V_{\delta}=V_0-\delta \phi.
$$
It suffices to show that $V_{\delta}\in C(q,r)$. Suppose that $q\in \partial \overline H_{V_{\delta}}(p')$ for some $p'\in \Rset^2$. The goal is to verify that
$$
\mathcal{M}_{p',V_\delta}\not=\Bbb T^2.
$$

\medskip

\nit {\bf Claim:} 
$$
\overline H_{V_\delta}(p_0)=\overline H_{V_0}(p_0).
$$
Since $V_\delta\leq V$, owing to the inf-max formula (\ref{inf-max}), 
$$
\overline H_{V_\delta}(p_0)\leq \overline H_{V_0}(p_0)
$$
So it suffices to show that 
$$
\overline H_{V_\delta}(p_0)\geq \overline H_{V_0}(p_0).
$$
Due to Corollary \ref{p-calib},
$$
\begin{array}{ll}
\int_{0}^{T_0}{1\over 2}|\dot \xi_\alpha|^2-V_{\delta}(\xi_\alpha(t))+\overline H_{V_\delta}(p_0)\,dt&\geq p_0\cdot (m,n)\\[5mm]
&=\int_{0}^{T_0}{1\over 2}|\dot \xi_\alpha|^2-V_0(\xi_\alpha(t))+\overline H_{V_0}(p_0)\,dt.
\end{array}
$$
Hence $\overline H_{V_\delta}(p)\geq \overline H_{V_0}(p)$. Therefore our claim holds and $\xi_\alpha$ is also a periodic orbit on $\mathcal{M}_{p_0,V_\delta}$ and all periodic orbits on $\mathcal{M}_{p_0,V_\delta}$ have the first homology class $(m,n)$. Similarly, $\xi_\beta$ is also a periodic orbit on $\mathcal{M}_{p_0,V_\delta}$. Hence
$$
q\in [\alpha q, \beta q]\subset \partial \overline H_{V_\delta}(p_0).
$$
Therefore $\overline H_{V_\delta}$ is linear on $[p_0,p']$. Owing to Property 3 in section 1 and (\ref{Matherequal}), 
$$
\mathcal{M}_{p_0,V_\delta}=\mathcal{M}_{p',V_\delta}.
$$
So we just need to prove that
$$
\mathcal{M}_{p_0,V_\delta}\not=\Bbb T^2.
$$
We argue by contradiction. Assume that 
$$
\mathcal{M}_{p_0,V_\delta}=\Bbb T^2.
$$
Let $\eta_0:[0,\tilde T_0]\to \Bbb T^2$ be the periodic orbit on $\mathcal{M}_{p_0,V_{\delta}}$ with $\eta_0(0)=x_0$. Here $\tilde T_0>0$ is the minimum periodic. Then
$$
\begin{array}{ll}
p_0\cdot (m,n)&=\int_{0}^{\tilde T}{1\over 2}|\dot \eta_0|^2-V_{\delta}(\eta_0)+\overline H_{V_\delta}(p_0)\,dt\\[5mm]
&=\int_{0}^{\tilde T}{1\over 2}|\dot \eta_0|^2-V_{\delta}(\eta_0)+\overline H_{V_0}(p_0)\,dt\\[5mm]
&>\int_{0}^{\tilde T}{1\over 2}|\dot \eta_0|^2-V_0(\eta_0)+\overline H_{V_0}(p_0)\,dt.
\end{array}
$$
This contradicts to Corollary \ref{p-calib}.\qed

\medskip

The residual set $\mathcal {G}$ in Theorem \ref{main2} is chosen as 
\be\label{residue}
\mathcal {G}=\cap_{q\in \Qset^2, \ r\in \Qset}C(q,r).
\ee
Here $\Qset^2$ is the collection of points on $\Rset^2$ whose both coordinates are rational numbers.

\medskip

\nit{\bf Acknowledgement} The author would like to thank Jinxin Xue for extremely helpful discussions about Aubry-Mather theory and numerous suggestions to improve the presentation of this paper. The author first thought about this problem when he was visiting Nanjing University in China in 2019. Fruitful discussions with Chong-qing Cheng and Wei Cheng at Nanjing University and their hospitality are greatly appreciated.  My gratitude also goes to Hongwei Gao for helping draw figures in this paper. The author also would like to thank Professor Lawrence Craig Evans for his encouragement to write down this result. 
\bibliographystyle{plain}

\begin{thebibliography}{99}



\bibitem{B1}
V. Bangert, 
{\em Mather Sets for Twist Maps and Geodesics on Tori}, Dynamics Reported, 1988, Volume 1.

\bibitem{B2}
V. Bangert, 
{\em Geodesic rays, Busemann functions and monotone twist maps}, Calculus of Variations and Partial Differential Equations
January 1994, Volume 2, Issue 1, pp 49--63. 


\bibitem{BC}

P. Bernard and G. Contreras, {\em A generic property of families of Lagrangian
systems}, Annals of Math., 167 (2008), 109--1108.


\bibitem{CZ}
C.Q. Cheng, M. Zhou, {\em Global normally hyperbolic invariant
cylinders in Lagrangian systems}, Math. Res. Lett, Volume 23, No. 3, 685--705, 2016.

\bibitem{C}

M. J. Carneiro, 
{\em On minimizing measures of the action of autonomous Lagrangians},
Nonlinearity 8 (1995) 1077--1085.




\bibitem{Co}
M. Concordel, {\em Periodic homogenisation of Hamilton-Jacobi equations II: eikonal equations}, Proc.
Roy. Soc. Edinburgh 127 (1997), 665--689.



\bibitem{W-E}
W. E,
\emph{Aubry-Mather theory and periodic solutions of the forced Burgers equation},
Comm. Pure Appl. Math. 52 (1999), no. 7, 811--828. 

\bibitem{G2002}

D. A. Gomes, \emph{Viscosity solutions of Hamilton-Jacobi equations, and asymptotics for Hamiltonian systems}, Calc. Var. 14, 345--357 (2002).

\bibitem{EG}
L. C. Evans, D. Gomes, 
\emph{Effective Hamiltonians and Averaging for
Hamiltonian Dynamics. I}, Arch. Ration. Mech. Anal. 157 (2001), no. 1, 1--33.


\bibitem{F}
A. Fathi, 
\emph{The weak KAM theorem in Lagrangian dynamics},
Cambridge University Press (2004).

\bibitem{H}
G. A. Hedlund, {\em Geodesics on a two-dimensional Riemannian manifold with periodic coefficients}, Ann. of Math. 33 (1932), 719--739.


\bibitem{JTY}
Wenjia Jing, Hung V. Tran, Yifeng Yu, {\em Effective fronts of polytope shapes}, to appear in Minimax Theory and its Applications.



\bibitem{Lions}

P. L. Lions, {\em Generalized solutions of Hamilton-Jacobi equations}, Research Notes in Mathematics, 69, 1982.

\bibitem{LPV} 
P. L. Lions, G. C. Papanicolaou, and S. R. S. Varadhan, 
{\em Homogenization of Hamilton-Jacobi equation}, Unpublished preprint, 1987. 


\bibitem{Ma1}

J. Mather, {\em Action minimizing invariant measures for positive definite
Lagrangian systems}, Math. Z., 207 (1991), 169--207.


\end{thebibliography}

\end{document}